\documentclass[11pt]{amsart}
\usepackage{amssymb}
\usepackage{amsmath}
\usepackage[active]{srcltx}
\usepackage{t1enc}
\usepackage[latin2]{inputenc}
\usepackage{verbatim}
\usepackage{amsmath,amsfonts,amssymb,amsthm}
\usepackage[mathcal]{eucal}
\usepackage{enumerate}
\usepackage[centertags]{amsmath}
\usepackage{graphics}

\newtheorem{theorem}{Theorem}

\newtheorem{lemma}{Lemma}

\begin{document}
\author{George Tephnadze}
\title[Fej\'er means ]{On The maximal operators of Vilenkin-Fejér means on
Hardy spaces}
\address{G. Tephnadze, Department of Mathematics, Faculty of Exact and
Natural Sciences, Tbilisi State University, Chavchavadze str. 1, Tbilisi
0128, Georgia}
\email{giorgitephnadze@gmail.com}
\date{}
\maketitle

\begin{abstract}
The main aim of this paper is to prove that when $0<p<1/2$ the maximal
operator $\overset{\sim }{\sigma }_{p}^{\ast }f:=\underset{n\in \mathbb{N}}{%
\sup }\frac{\left\vert \sigma _{n}f\right\vert }{\left( n+1\right) ^{1/p-2}}$
is bounded from the martingale Hardy space $H_{p}$ to the space $L_{p},$
where $\sigma _{n}$ is $n$-th Fejér mean with respect to bounded Vilenkin
system.
\end{abstract}

\date{}

\textbf{2000 Mathematics Subject Classification.} 42C10.

\textbf{Key words and phrases:} Vilenkin system, Fejér means, martingale
Hardy space.

\section{INTRODUCTION}

In one-dimensional case the weak type inequality for maximal operator of Fejé%
r means for trigonometric system can be found in Zygmund \cite{Zy}, in
Schipp \cite{Sc} for Walsh system and in Pál, Simon \cite{PS} for bounded
Vilenkin system. Fujii \cite{Fu} and Simon \cite{Si2} verified that the $%
\sigma ^{\ast ,w}$ is bounded from $H_{1}$ to $L_{1},$ where $\sigma ^{\ast
,w}$ denotes the maximal operator of Fejér means of Walsh-Fourier series.
Weisz \cite{We2} generalized this result and proved the boundedness of $%
\sigma ^{\ast ,w}$ from the martingale Hardy space $H_{p}$ to the space $%
L_{p},$ for $1/2<p\leq 1$. Simon \cite{Si1} gave a counterexample, which
shows that boundedness of $\sigma ^{\ast ,w}$ does not hold for $0<p<1/2.$
The counterexample for $\sigma ^{\ast ,w}$ when $p=1/2$ is due to Goginava
\cite{GoPubl} (see also \cite{BGG2,tep1}). In the endpoint case $p=1/2$ two
positive results were showed. Weisz \cite{we4} proved that $\sigma ^{\ast
,w} $ is bounded from the Hardy space $H_{1/2}$ to the space $L_{1/2,\infty }
$. Goginava \cite{GoSzeged} proved that the maximal operator $\widetilde{%
\sigma }$ $^{\ast ,w\,}$defined by
\begin{equation*}
\widetilde{\sigma }^{\ast ,w}f:=\sup_{n\in \mathbb{N}}\frac{\left\vert
\sigma _{n}^{w}f\right\vert }{\log ^{2}\left( n+1\right) }
\end{equation*}%
is bounded from the Hardy space $H_{1/2}$ to the space $L_{1/2},$ where $%
\sigma _{n}^{w}$ is $n$-th Fejér means of Walsh-Fourier series. He also
proved, that for any nondecreasing function $\varphi :\mathbb{N}\rightarrow
\lbrack 1,$ $\infty ),$ satisfying the condition
\begin{equation*}
\overline{\lim_{n\rightarrow \infty }}\frac{\log ^{2}\left( n+1\right) }{%
\varphi \left( n\right) }=+\infty ,
\end{equation*}%
the maximal operator
\begin{equation*}
\sup_{n\in \mathbb{N}}\frac{\left\vert \sigma _{n}^{w}f\right\vert }{\varphi
\left( n\right) }
\end{equation*}%
is not bounded from the Hardy space $H_{1/2}$ to the space $L_{1/2}.$

For Walsh-Kaczmarz system analogical theorem was proved in \cite{GNCz} and
for bounded Vilenkin system in \cite{tep2}.

The main aim of this paper is to prove that when $0<p<1/2$ the maximal
operator

\begin{equation}
\widetilde{\sigma }_{p}^{\ast }f\,:=\underset{n\in \mathbb{N}}{\sup }\frac{%
\left\vert \sigma _{n}f\right\vert }{\left( n+1\right) ^{1/p-2}}
\label{cond}
\end{equation}%
is bounded from the Hardy space $H_{p}$ to the space $L_{p}$ (see Theorem
1), where $\sigma _{n}$ is Fejér means of bounded Vilenkin-Fourier series.

We also prove that for any nondecreasing function $\varphi :\mathbb{N}%
\rightarrow \lbrack 1,$ $\infty ),$ satisfying the condition
\begin{equation}
\overline{\lim_{n\rightarrow \infty }}\frac{\left( n+1\right) ^{1/p-2}}{%
\varphi \left( n\right) }=+\infty ,  \label{cond2}
\end{equation}%
the maximal operator
\begin{equation*}
\sup_{n\in \mathbb{N}}\frac{\left\vert \sigma _{n}f\right\vert }{\varphi
\left( n\right) }
\end{equation*}%
is not bounded from the Hardy space $H_{p}$ to the space $L_{p,\infty }$
when $0<p<1/2.$ Actually, we prove a stronger result (see Theorem 2) than
the unboundedness of the maximal operator $\widetilde{\sigma }_{p}^{\ast }$
from the Hardy space $H_{p}$ to the spaces $L_{p,\infty }.$ In particular,
we prove that under condition (\ref{cond2}) there exists a martingale $f\in
H_{p}$ $\left( 0<p<1/2\right) $ such that

\begin{equation*}
\underset{n\in \mathbb{N}}{\sup }\left\Vert \frac{\sigma _{n}f}{\varphi
\left( n\right) }\right\Vert _{L_{p,\infty }}=\infty .
\end{equation*}

\section{Definitions and Notations}

Let $\mathbb{N}_{+}$ denote the set of the positive integers, $\mathbb{N}:=%
\mathbb{N}_{+}\cup \{0\}.$

Let $m:=(m_{0,}...,$ $m_{n},...)$ denote a sequence of the positive integers
not less than 2.

Denote by
\begin{equation*}
Z_{m_{k}}:=\{0,1,...,\text{ }m_{k}-1\}
\end{equation*}
the additive group of integers modulo $m_{k}.$

Define the group $G_{m},$ as the complete direct product of the group $%
Z_{m_{j}},$ with the product of the discrete topologies of $Z_{m_{j}}$ $^{,}$%
s.

The direct product $\mu $ of the measures
\begin{equation*}
\mu _{k}\left( \{j\}\right) :=1/m_{k},\text{ \qquad }(j\in Z_{m_{k}}),
\end{equation*}
is the Haar measure on $G_{m},$ with $\mu \left( G_{m}\right) =1.$

If $\sup_{n}m_{n}<\infty $, then we call $G_{m}$ a bounded Vilenkin group.
If the generating sequence $m$ is not bounded, then $G_{m}$ is said to be an
unbounded Vilenkin group. \textbf{In this paper we discuss bounded Vilenkin
groups only.}

The elements of $G_{m}$ represented by sequences
\begin{equation*}
x:=(x_{0},x_{1},...,x_{j},...),\qquad \left( \text{ }x_{k}\in
Z_{m_{k}}\right) .
\end{equation*}

It is easy to give a base, for the neighborhood of $\ x\in G_{m}:$
\begin{equation*}
I_{0}\left( x\right) :=G_{m},
\end{equation*}%
\begin{equation*}
I_{n}(x):=\{y\in G_{m}\mid y_{0}=x_{0},...,\text{ }y_{n-1}=x_{n-1}\},\text{ }%
(n\in \mathbb{N}).
\end{equation*}%
Denote $I_{n}:=I_{n}\left( 0\right) ,$ for $n\in \mathbb{N}$ and $\overset{-}%
{I_{n}}:=G_{m}$ $\backslash $ $I_{n}.$

Let

\begin{equation*}
e_{n}:=\left( 0,...,0,x_{n}=1,0,...\right) \in G_{m},\qquad \left( n\in
\mathbb{N}\right) .
\end{equation*}

Denote
\begin{equation*}
I_{N}^{k,l}:=\left\{
\begin{array}{l}
\text{ }I_{N}(0,...,0,x_{k}\neq 0,0,...,0,x_{l}\neq 0,x_{l+1\text{ }%
},...,x_{N-1\text{ }},...),\text{ } \\
\text{ where }x_{i\text{ }}\in Z_{m_{i}},\text{ }i\geq l+1,\text{ for }k<l<N,
\\
\text{ }I_{N}(0,...,0,x_{k}\neq 0,0,...,,x_{N-1\text{ }}=0,\text{ }x_{N\text{
}},...),\text{ } \\
\text{where }x_{i\text{ }}\in Z_{m_{i}},\text{ }i\geq N,\text{ for }l=N.%
\end{array}%
\text{ }\right.
\end{equation*}%
and
\begin{equation*}
I_{N}^{k,\alpha ,l,\beta }:=I_{N}(0,...,0,x_{k}=\alpha ,0,...,0,x_{l}=\beta
,x_{l+1,...,\text{ }}x_{N-1\text{ }}),\text{ \ }k<l<N\text{\ },
\end{equation*}%
where $x_{i\text{ }}\in Z_{m_{i}},$ $i\geq l+1.$

It is evident
\begin{equation}
I_{N}^{k,l}=\overset{m_{k}-1}{\underset{\alpha =1}{\bigcup }}\overset{m_{l}-1%
}{\underset{\beta =1}{\bigcup }}I_{N}^{k,\alpha ,l,\beta }  \label{1}
\end{equation}%
and
\begin{equation}
\overset{-}{I_{N}}=\left( \overset{N-2}{\underset{k=0}{\bigcup }}\overset{N-1%
}{\underset{l=k+1}{\bigcup }}I_{N}^{k,l}\right) \bigcup \left( \underset{k=1}%
{\bigcup\limits^{N-1}}I_{N}^{k,N}\right) .  \label{2}
\end{equation}

If we define the so-called generalized number system, based on $m$ in the
following way :
\begin{equation*}
M_{0}:=1,\text{ \qquad }M_{k+1}:=m_{k}M_{k\text{ }},\ \qquad (k\in \mathbb{N}%
)
\end{equation*}%
then every $n\in \mathbb{N}$ can be uniquely expressed as $n=\overset{\infty
}{\underset{k=0}{\sum }}n_{j}M_{j},$ where $n_{j}\in Z_{m_{j}}$ $~(j\in
\mathbb{N})$ and only a finite number of $n_{j}`$s differ from zero. Let $%
\left\vert n\right\vert :=\max $ $\{j\in \mathbb{N},$ $n_{j}\neq 0\}.$

It is easy to show that%
\begin{equation}
\overset{l}{\underset{A=0}{\sum }}M_{A}\leq cM_{l}.  \label{3a}
\end{equation}

Denote by $L_{1}\left( G_{m}\right) $ the usual (one dimensional) Lebesque
space.

Next, we introduce on $G_{m}$ an ortonormal system which is called the
Vilenkin system.

At first define the complex valued function $r_{k}\left( x\right)
:G_{m}\rightarrow C,$ the generalized Rademacher functions as
\begin{equation*}
r_{k}\left( x\right) :=\exp \left( 2\pi ix_{k}/m_{k}\right) ,\text{ \qquad }%
\left( i^{2}=-1,\text{ }x\in G_{m},\text{ }k\in \mathbb{N}\right) .
\end{equation*}

Now define the Vilenkin system $\psi :=(\psi _{n}:n\in \mathbb{N})$ on $%
G_{m} $ as:
\begin{equation*}
\psi _{n}(x):=\overset{\infty }{\underset{k=0}{\Pi }}r_{k}^{n_{k}}\left(
x\right) ,\text{ \qquad }\left( n\in \mathbb{N}\right) .
\end{equation*}

Specifically, we call this system the Walsh one if $m\equiv 2$.

The Vilenkin system is ortonormal and complete in $L_{2}\left( G_{m}\right)
\,$\cite{AVD,Vi}.

Now we introduce analogues of the usual definitions in Fourier-analysis.

If $f\in L_{1}\left( G_{m}\right) $ we can establish the the Fourier
coefficients, the partial sums of the Fourier series, the Fejér means, the
Dirichlet and Fejér kernels with respect to the Vilenkin system $\psi $ in
the usual manner:
\begin{eqnarray*}
\widehat{f}\left( n\right) &:&=\int_{G_{m}}f\overline{\psi }_{n}d\mu ,\text{%
\thinspace \qquad }\left( \text{ }n\in \mathbb{N}\text{ }\right) , \\
S_{n}f &:&=\sum_{k=0}^{n-1}\widehat{f}\left( k\right) \psi _{k}\ ,\text{
\qquad }\left( \text{ }n\in \mathbb{N}_{+},\text{ }S_{0}f:=0\right) , \\
\sigma _{n}f &:&=\frac{1}{n}\sum_{k=0}^{n-1}S_{k}f,\text{ \qquad }\left(
\text{ }n\in \mathbb{N}_{+}\text{ }\right) , \\
D_{n} &:&=\sum_{k=0}^{n-1}\psi _{k\text{ }},\text{ \qquad }\left( \text{ }%
n\in \mathbb{N}_{+}\text{ }\right) , \\
K_{n} &:&=\frac{1}{n}\overset{n-1}{\underset{k=0}{\sum }}D_{k},\text{%
\thinspace \qquad }\left( \text{ }n\in \mathbb{N}_{+}\text{ }\right) .
\end{eqnarray*}

Recall that
\begin{equation}
\quad \hspace*{0in}D_{M_{n}}\left( x\right) =\left\{
\begin{array}{l}
\text{ }M_{n},\,\text{\thinspace \thinspace if\thinspace\ \ \thinspace\ }%
x\in I_{n}, \\
\text{ ~}0\text{ },\text{ \thinspace \thinspace \thinspace if \thinspace\ \
\thinspace }x\notin I_{n}.%
\end{array}%
\right.  \label{3}
\end{equation}

It is well-known that
\begin{equation}
\sup_{n}\int_{G_{m}}\left| K_{n}\left( x\right) \right| d\mu \left( x\right)
\leq c<\infty ,  \label{4}
\end{equation}
and
\begin{equation}
n\left| K_{n}\left( x\right) \right| \leq c\sum_{A=0}^{\left| n\right|
}M_{A}\left| K_{M_{A}}\left( x\right) \right| .  \label{5}
\end{equation}
\vspace{0pt}

The norm (or quasinorm) of the space $L_{p}(G_{m})$ is defined by \qquad
\qquad \thinspace\
\begin{equation*}
\left\Vert f\right\Vert _{p}:=\left( \int_{G_{m}}\left\vert f(x)\right\vert
^{p}d\mu (x)\right) ^{1/p},\qquad \left( 0<p<\infty \right) .
\end{equation*}%
The space $L_{p,\infty }\left( G_{m}\right) $ consists of all measurable
functions $f,$ for which

\begin{equation*}
\left\Vert f\right\Vert _{L_{p,\infty }(G_{m})}:=\underset{\lambda >0}{\sup }%
\lambda ^{p}\mu \left( f>\lambda \right) <+\infty .
\end{equation*}

The $\sigma -$algebra generated by the intervals $\left\{ I_{n}\left(
x\right) :x\in G_{m}\right\} $ will be denoted by $\digamma _{n}$ $\left(
n\in \mathbb{N}\right) .$ Denote by $f=\left( f^{\left( n\right) },n\in
\mathbb{N}\right) $ a martingale with respect to $\digamma _{n}$ $\left(
n\in \mathbb{N}\right) $ (for details see e.g. \cite{We1}). The maximal
function of a martingale $f$ is defend by \qquad
\begin{equation*}
f^{\ast }=\sup_{n\in \mathbb{N}}\left\vert f^{\left( n\right) }\right\vert .
\end{equation*}

In case $f\in L_{1},$ the maximal functions are also be given by
\begin{equation*}
f^{\ast }\left( x\right) =\sup_{n\in \mathbb{N}}\frac{1}{\left\vert
I_{n}\left( x\right) \right\vert }\left\vert \int_{I_{n}\left( x\right)
}f\left( u\right) \mu \left( u\right) \right\vert .
\end{equation*}

For $0<p<\infty $ the Hardy martingale spaces $H_{p}$ $\left( G_{m}\right) $
consist of all martingales for which
\begin{equation*}
\left\Vert f\right\Vert _{H_{p}}:=\left\Vert f^{\ast }\right\Vert
_{p}<\infty .
\end{equation*}

If $f\in L_{1},$ then it is easy to show that the sequence $\left(
S_{M_{n}}\left( f\right) :n\in \mathbb{N}\right) $ is a martingale. If $%
f=\left( f^{\left( n\right) },n\in \mathbb{N}\right) $ is a martingale, then
the Vilenkin-Fourier coefficients must be defined in a slightly different
manner: $\qquad \qquad $
\begin{equation*}
\widehat{f}\left( i\right) :=\lim_{k\rightarrow \infty
}\int_{G_{m}}f^{\left( k\right) }\left( x\right) \overline{\Psi }_{i}\left(
x\right) d\mu \left( x\right) .
\end{equation*}%
\qquad \qquad \qquad \qquad

The Vilenkin-Fourier coefficients of $f\in L_{1}\left( G_{m}\right) $ are
the same as those of the martingale $\left( S_{M_{n}}\left( f\right) :n\in
\mathbb{N}\right) $ obtained from $f$ .

For the martingale $f$ we consider maximal operators
\begin{eqnarray*}
\sigma ^{\ast }f &:&=\sup_{n\in \mathbb{N}}\left\vert \sigma
_{n}f\right\vert ,\text{ } \\
\widetilde{\sigma }^{\ast }f &:&=\sup_{n\in \mathbb{N}}\frac{\left\vert
\sigma _{n}f\right\vert }{\log ^{2}\left( n+1\right) },\text{ } \\
\widetilde{\sigma }_{p}^{\ast }f &:&=\sup_{n\in \mathbb{N}}\frac{\left\vert
\sigma _{n}f\right\vert }{\left( n+1\right) ^{1/p-2}}.
\end{eqnarray*}

A bounded measurable function $a$ is p-atom, if there exist interval $I$,
such that \qquad
\begin{eqnarray*}
a)\text{ }\int_{I}ad\mu &=&0, \\
b)\ \left\Vert a\right\Vert _{\infty } &\leq &\mu \left( I\right) ^{-1/p}, \\
c)\text{ supp}\left( a\right) &\subset &I.
\end{eqnarray*}

\section{Formulation of Main Results}

\begin{theorem}
Let $0<p<1/2.$ Then the \bigskip maximal operator $\widetilde{\sigma }%
_{p}^{*}$ is bounded from the Hardy martingale space $H_{p}\left(
G_{m}\right) $ to the space $L_{p}\left( G_{m}\right) .$
\end{theorem}

\begin{theorem}
Let $\varphi :\mathbb{N}\rightarrow \lbrack 1,$ $\infty )$ be a
nondecreasing function, satisfying the condition
\end{theorem}

\begin{equation}
\overline{\lim_{n\rightarrow \infty }}\frac{\left( n+1\right) ^{1/p-2}}{%
\varphi \left( n\right) }=+\infty ,  \label{6}
\end{equation}%
then
\begin{equation*}
\sup_{n\in \mathbb{N}}\left\Vert \frac{\sigma _{n}f}{\varphi \left( n\right)
}\right\Vert _{L_{p,\infty }}=\infty .
\end{equation*}

\section{AUXILIARY PROPOSITIONS}

\begin{lemma}
\cite{We3} Suppose that an operator $T$ is sublinear and for some $0<p\leq 1$
\end{lemma}

\begin{equation*}
\int\limits_{\overset{-}{I}}\left| Ta\right| ^{p}d\mu \leq c_{p}<\infty ,
\end{equation*}
for every $p$-atom $a$, where $I$ denote the support of the atom. If $T$ is
bounded from $L_{\infty \text{ }}$ to $L_{\infty },$ then
\begin{equation*}
\left\| Tf\right\| _{L_{p}\left( G_{m}\right) }\leq c_{p}\left\| f\right\|
_{H_{p}\left( G_{m}\right) }.
\end{equation*}
\bigskip

\begin{lemma}
\cite{BGG,GoAMH} Let $2<A\in \mathbb{N}_{+},$ $k\leq s<A$ and $%
q_{A}=M_{2A}+M_{2A-2}+...+M_{2}+M_{0}.$ Then
\begin{equation*}
q_{A-1}\left\vert K_{q_{A-1}}(x)\right\vert \geq \frac{M_{2k}M_{2s}}{4},
\end{equation*}%
for
\begin{eqnarray*}
x &\in &I_{2A}\left( 0,...,x_{2k}\neq 0,0,...,0,x_{2s}\neq
0,x_{2s+1},...x_{2A-1}\right) , \\
k &=&0,\text{ }1,...,\text{ }A-3.\qquad s=k+2,\text{ }k+3,...,\text{ }A-1.
\end{eqnarray*}
\end{lemma}

\begin{lemma}
\cite{gat} Let $A>t,$ $t,A\in \mathbb{N},$ $z\in I_{t}\backslash $ $I_{t+1}$%
. Then
\end{lemma}

$\quad \hspace*{0in}$
\begin{equation*}
K_{M_{A}}\left( z\right) =\left\{
\begin{array}{c}
\text{ }0,\text{\qquad if }z-z_{t}e_{t}\notin I_{A}, \\
\text{ }\frac{M_{t}}{1-r_{t}\left( z\right) },\text{\qquad if }%
z-z_{t}e_{t}\in I_{A}.%
\end{array}
\right.
\end{equation*}

\begin{lemma}
Let $x\in I_{N}^{k,l}$ $,$\qquad $k=0,...,N-1,$ $l=k+1,...,N.$ Then
\end{lemma}

\begin{equation*}
\int_{I_{N}}\left| K_{n}\left( x-t\right) \right| d\mu \left( t\right) \leq
\frac{cM_{l}M_{k}}{M_{N}^{2}},\,\,\,\,\text{when\thinspace \thinspace
\thinspace }n\geq M_{N}.
\end{equation*}
\textbf{Proof.} Let $x\in I_{N}^{k,\alpha ,l,\beta }$. Then applying Lemma 3
we have
\begin{equation*}
K_{M_{A}}\left( x\right) =0,\,\,\text{when \thinspace \thinspace }A>l.
\end{equation*}

Let $k<A\leq l$. Then we get
\begin{equation}
\left| K_{M_{A}}\left( x\right) \right| =\frac{M_{k}}{\left| 1-\text{ }%
r_{k}\left( x\right) \right| }\leq \frac{m_{k}M_{k}}{2\pi \text{ }\alpha }.
\label{7}
\end{equation}

Let $x\in I_{N}^{k,l},$ for $0\leq k<l\leq N-1$ and $t\in I_{N}.$ Since $%
x-t\in $ $I_{N}^{k,l}$ and $n\geq M_{N},$ combining (\ref{1}), (\ref{3a}), (%
\ref{5}) and (\ref{7}) we obtain

\begin{equation*}
n\left\vert K_{n}\left( x\right) \right\vert \leq c\overset{l}{\underset{A=0}%
{\sum }}M_{A}M_{k}\leq cM_{k}M_{l}
\end{equation*}%
and

\begin{equation}
\int_{I_{N}}\left| K_{n}\left( x-t\right) \right| d\mu \left( t\right) \leq
\frac{cM_{k}M_{l}}{M_{N}^{2}}.  \label{8}
\end{equation}

Let $x\in I_{N}^{k,N}$ , then applying (\ref{5}) we have

\begin{eqnarray}
&&\int_{I_{N}}n\left| K_{n}\left( x-t\right) \right| d\mu \left( t\right)
\label{9} \\
&\leq &\underset{A=0}{\overset{\left| n\right| }{\sum }}M_{A}\int_{I_{N}}%
\left| K_{M_{A}}\left( x-t\right) \right| d\mu \left( t\right) .  \notag
\end{eqnarray}

Let

\begin{equation*}
\left\{
\begin{array}{l}
x=\left( 0,...,0,x_{k}\neq 0,0,...0,x_{N},x_{N+1},x_{q},...,x_{\left|
n\right| -1},...\right) , \\
t=\left( 0,...,0,x_{N},...x_{q-1},t_{q}\neq x_{q},t_{q+1},...,t_{\left|
n\right| -1},...\right) ,\,\,q=N,...,\left| n\right| -1.\text{ }%
\end{array}
\right.
\end{equation*}

Using Lemma 3 in (\ref{9}) it is easy to show that

\begin{eqnarray}
&&\int_{I_{N}}\left| K_{n}\left( x-t\right) \right| d\mu \left( t\right)
\label{10} \\
&\leq &\frac{c}{n}\underset{A=0}{\overset{q-1}{\sum }}M_{A}%
\int_{I_{N}}M_{k}d\mu \left( t\right) \leq \frac{cM_{k}M_{q}}{nM_{N}}\leq
\frac{cM_{k}}{M_{N}}.  \notag
\end{eqnarray}

Let

\begin{equation*}
\left\{
\begin{array}{l}
\text{ }x=\left( 0,...,0,x_{m}\neq
0,0,...,0,x_{N},x_{N+1},x_{q},...,x_{\left\vert n\right\vert -1},...\right) ,%
\text{ } \\
t=\left( 0,0,...,x_{N},...,x_{_{\left\vert n\right\vert -1}},...\right) .%
\end{array}%
\right.
\end{equation*}

If we apply Lemma 3 in (\ref{9}) we obtain
\begin{eqnarray}
&&\int_{I_{N}}\left| K_{n}\left( x-t\right) \right| d\mu \left( t\right)
\label{11} \\
&\leq &\frac{c}{n}\overset{\left| n\right| -1}{\underset{A=0}{\sum }}%
M_{A}\int_{I_{N}}M_{k}d\mu \left( t\right) \leq \frac{cM_{k}}{M_{N}}.  \notag
\end{eqnarray}

Combining (\ref{8}), (\ref{10}) and (\ref{11}) we complete the proof of \
Lemma 4.

\section{Proofs of the Theorems}

\textbf{Proof of Theorem 1. }By Lemma 1, the proof of Theorem 1 will be
complete, if we show that

\begin{equation*}
\int\limits_{\overline{I}_{N}}\left( \underset{n\in \mathbb{N}}{\sup }\frac{%
\left\vert \sigma _{n}a\right\vert }{\left( n+1\right) ^{1/p-2}}\right)
^{p}d\mu \leq c<\infty ,
\end{equation*}%
for every p-atom $a,$ where $I$ denotes the support of the atom$.$ The
boundedness of $\ \sup_{n\in \mathbb{N}}\left\vert \sigma _{n}f\right\vert
/\left( n+1\right) ^{1/p-2}$ from $L_{\infty }$ to $L_{\infty }$ follows
from (\ref{4}).

Let $a$ be an arbitrary p-atom, with support$\ I$ and $\mu \left( I\right)
=M_{N}^{-1}.$ We may assume that $I=I_{N}.$ It is easy to see that $\sigma
_{n}\left( a\right) =0,$ when $n\leq M_{N}$. Therefore, we can suppose that $%
n>M_{N}$.

Since $\left\| a\right\| _{\infty }\leq cM_{N}^{1/p}$ we can write
\begin{eqnarray*}
&&\frac{\left| \sigma _{n}\left( a\right) \right| }{\left( n+1\right)
^{1/p-2}} \\
&\leq &\frac{1}{\left( n+1\right) ^{1/p-2}}\int_{I_{N}}\left| a\left(
t\right) \right| \left| K_{n}\left( x-t\right) \right| d\mu \left( t\right)
\\
&\leq &\frac{\left\| a\right\| _{\infty }}{\left( n+1\right) ^{1/p-2}}%
\int_{I_{N}}\left| K_{n}\left( x-t\right) \right| d\mu \left( t\right) \\
&\leq &\frac{cM_{N}^{1/p}}{\left( n+1\right) ^{1/p-2}}\int_{I_{N}}\left|
K_{n}\left( x-t\right) \right| d\mu \left( t\right) .
\end{eqnarray*}

Let $x\in I_{N}^{k,l},\,0\leq k<l\leq N.$ From Lemma 4 we get
\begin{equation}
\frac{\left\vert \sigma _{n}\left( a\right) \right\vert }{\left( n+1\right)
^{1/p-2}}\leq \frac{cM_{N}^{1/p}}{M_{N}^{1/p-2}}\frac{M_{l}M_{k}}{M_{N}^{2}}%
=cM_{l}M_{k}.  \label{12}
\end{equation}%
Combining (\ref{2}) and (\ref{12}) we obtain
\begin{eqnarray*}
&&\int_{\overline{I_{N}}}\left\vert \sigma ^{\ast }a\left( x\right)
\right\vert ^{p}d\mu \left( x\right) \\
&=&\overset{N-2}{\underset{k=0}{\sum }}\overset{N-1}{\underset{l=k+1}{\sum }}%
\sum\limits_{x_{j}=0,j\in
\{l+1,...,N-1\}}^{m_{j-1}}\int_{I_{N}^{k,l}}\left\vert \sigma ^{\ast
}a\left( x\right) \right\vert ^{p}d\mu \left( x\right) \\
&&+\overset{N-1}{\underset{k=0}{\sum }}\int_{I_{N}^{k,N}}\left\vert \sigma
^{\ast }a\left( x\right) \right\vert ^{p}d\mu \left( x\right) \\
&\leq &c\overset{N-2}{\underset{k=0}{\sum }}\overset{N-1}{\underset{l=k+1}{%
\sum }}\frac{m_{l+1}...m_{N-1}}{M_{N}}\left( M_{l}M_{k}\right) ^{p} \\
&&+c\overset{N-1}{\underset{k=0}{\sum }}\frac{1}{M_{N}}\left(
M_{N}M_{k}\right) ^{p} \\
&\leq &c\overset{N-2}{\underset{k=0}{\sum }}\overset{N-1}{\underset{l=k+1}{%
\sum }}\frac{\left( M_{l}M_{k}\right) ^{p}}{M_{l}}+c\overset{N-1}{\underset{%
k=0}{\sum }}\frac{M_{k}^{p}}{M_{N}^{1-p}}=I+II.
\end{eqnarray*}

Then

\begin{eqnarray*}
I &=&c\overset{N-2}{\underset{k=0}{\sum }}\overset{N-1}{\underset{l=k+1}{%
\sum }}\frac{1}{M_{l}^{1-2p}}\frac{\left( M_{l}M_{k}\right) ^{p}}{M_{l}^{2p}}
\\
&\leq &c\overset{N-2}{\underset{k=0}{\sum }}\overset{N-1}{\underset{l=k+1}{%
\sum }}\frac{1}{M_{l}^{1-2p}} \\
&\leq &c\overset{N-2}{\underset{k=0}{\sum }}\overset{N-1}{\underset{l=k+1}{%
\sum }}\frac{1}{2^{l\left( 1-2p\right) }} \\
&\leq &c\overset{N-2}{\underset{k=0}{\sum }}\frac{1}{2^{k\left( 1-2p\right) }%
}<c<\infty .
\end{eqnarray*}

It is evident

\begin{equation*}
II\leq \frac{c}{M_{N}^{1-2p}} <c<\infty.
\end{equation*}

Which complete the proof of Theorem 1.

\textbf{Proof of Theorem 2.} Let $0<p<1/2$ and$\ \left\{ \lambda _{k};\text{
}k\in \mathbb{N}_{+}\right\} $ be an increasing sequence of the positive
integers, such that
\begin{equation*}
\lim_{k\rightarrow \infty }\frac{\lambda _{k}^{1/p-2}}{\varphi \left(
\lambda _{k}\right) }=\infty .
\end{equation*}%
It is evident that for every $\lambda _{k},$ there exists a positive
integers $m_{k}^{,},$ such that $q_{_{m_{k}^{^{\prime }}}}<\lambda
_{k}<cq_{_{m_{k}^{,}}}.$ Since $\varphi \left( n\right) $ is nondecreasing
function, we have
\begin{eqnarray}
&&\overline{\underset{k\rightarrow \infty }{\lim }}\frac{%
M_{2m_{k}^{,}}^{1/p-2}}{\varphi \left( q_{m_{k}^{,}}\right) }  \label{13} \\
&\geq &c\overline{\underset{k\rightarrow \infty }{\lim }}\frac{%
q_{m_{k}^{,}}^{1/p-2}}{\varphi \left( q_{m_{k}^{,}}\right) }  \notag \\
&\geq &c\lim_{k\rightarrow \infty }\frac{\lambda _{k}^{1/p-2}}{\varphi
\left( \lambda _{k}\right) }=\infty .  \notag
\end{eqnarray}

Let$\ \left\{ n_{k};\text{ }k\in \mathbb{N}_{+}\right\} \subset \left\{
m_{k}^{,};\text{ }k\in \mathbb{N}_{+}\right\} $ such that
\begin{equation*}
\lim_{k\rightarrow \infty }\frac{M_{2n_{k}}^{1/p-2}}{\varphi \left(
q_{n_{k}}\right) }=\infty
\end{equation*}%
and

\begin{equation*}
f_{n_{k}}\left( x\right) =D_{M_{2n_{k}+1}}\left( x\right)
-D_{M_{_{2n_{k}}}}\left( x\right) ,\text{ \qquad }n_{k}\geq 3.
\end{equation*}

It is evident
\begin{equation*}
\widehat{f}_{n_{k}}\left( i\right) =\left\{
\begin{array}{l}
\text{ }1,\text{ if }i=M_{_{2n_{k}}},...,M_{2n_{k}+1}-1, \\
\text{ }0,\text{otherwise}.%
\end{array}%
\right.
\end{equation*}%
Then we can write
\begin{equation}
S_{i}f_{n_{k}}\left( x\right) =\left\{
\begin{array}{l}
D_{i}\left( x\right) -D_{M_{_{2n_{k}}}}\left( x\right) ,\text{ if }%
i=M_{_{2n_{k}}},...,M_{2n_{k}+1}-1, \\
\text{ }f_{n_{k}}\left( x\right) ,\text{ if }i\geq M_{2n_{k}+1}, \\
0,\text{ \qquad otherwise}.%
\end{array}%
\right.  \label{14}
\end{equation}

From (\ref{3}) we get
\begin{eqnarray}
&&\left\Vert f_{n_{k}}\right\Vert _{H_{p}}  \label{15} \\
&=&\left\Vert \sup\limits_{n\in \mathbb{N}}S_{M_{n}}\left( f_{n_{k}}\right)
\right\Vert _{L_{p}}  \notag \\
&=&\left\Vert D_{M_{2n_{k}+1}}-D_{M_{_{2n_{k}}}}\right\Vert _{L_{p}}  \notag
\\
&=&\left( \int_{I_{_{2n_{k}}}\backslash \text{ }%
I_{_{2n_{k}+1}}}M_{_{2n_{k}}}^{p}d\mu \left( x\right)
+\int_{I_{_{2n_{k}+1}}}\left( M_{_{2n_{k}+1}}-M_{_{2n_{k}}}\right) ^{p}d\mu
\left( x\right) \right) ^{1/p}  \notag \\
&=&\left( \frac{m_{_{2n_{k}}}-1}{M_{2n_{k}+1}}M_{_{2n_{k}}}^{p}+\frac{\left(
m_{_{2n_{k}}}-1\right) ^{p}}{M_{_{2n_{k}}+1}}M_{_{2n_{k}}}^{p}\right) ^{1/p}
\notag \\
&\leq &M_{_{2n_{k}}}^{1-1/p}.  \notag
\end{eqnarray}%
By (\ref{14}) we can write:
\begin{eqnarray*}
&&\frac{\left\vert \sigma _{q_{n_{k}}}f_{n_{k}}\left( x\right) \right\vert }{%
\varphi \left( q_{n_{k}}\right) } \\
&=&\frac{1}{\varphi \left( q_{n_{k}}\right) q_{n_{k}}}\left\vert \overset{%
q_{n_{k}}-1}{\underset{j=0}{\sum }}S_{j}f_{n_{k}}\left( x\right) \right\vert
\\
&=&\frac{1}{\varphi \left( q_{n_{k}}\right) q_{n_{k}}}\left\vert \overset{%
q_{n_{k}}-1}{\underset{j=M_{_{2n_{k}}}}{\sum }}S_{j}f_{n_{k}}\left( x\right)
\right\vert \\
&=&\frac{1}{\varphi \left( q_{n_{k}}\right) q_{n_{k}}}\left\vert \overset{%
q_{n_{k}}-1}{\underset{j=M_{_{2n_{k}}}}{\sum }}\left( D_{j}\left( x\right)
-D_{M_{_{2n_{k}}}}\left( x\right) \right) \right\vert \\
&=&\frac{1}{\varphi \left( q_{n_{k}}\right) q_{n_{k}}}\left\vert \overset{%
q_{n_{k}-1}-1}{\underset{j=0}{\sum }}\left( D_{j+M_{_{2n_{k}}}}\left(
x\right) -D_{M_{_{2n_{k}}}}\left( x\right) \right) \right\vert .
\end{eqnarray*}%
Since

\begin{equation*}
D_{_{j+M_{_{2n_{k}}}}}\left( x\right) -D_{M_{_{2n_{k}}}}\left( x\right)
=\psi _{M_{_{2n_{k}}}}D_{j},\text{ \qquad }\,j=1,2,..,M_{_{2n_{k}}}-1,
\end{equation*}%
we obtain
\begin{eqnarray*}
&&\frac{\left\vert \sigma _{q_{n_{k}}}f_{n_{k}}\left( x\right) \right\vert }{%
\varphi \left( q_{n_{k}}\right) } \\
&=&\frac{1}{\varphi \left( q_{n_{k}}\right) q_{n_{k}}}\left\vert \overset{%
q_{n_{k}-1}-1}{\underset{j=0}{\sum }}D_{j}\left( x\right) \right\vert \\
&=&\frac{1}{\varphi \left( q_{n_{k}}\right) }\frac{q_{n_{k}-1}}{q_{n_{k}}}%
\left\vert K_{q_{n_{k}-1}}\left( x\right) \right\vert .
\end{eqnarray*}

Let $x\in $ $I_{_{2n_{k}}}^{2s,2l}$. Using Lemma 2 we obtain
\begin{equation*}
\frac{\left\vert \sigma _{q_{n_{k}}}f_{n_{k}}\left( x\right) \right\vert }{%
\varphi \left( q_{n_{k}}\right) }\geq \frac{cM_{2s}M_{2l}}{%
M_{_{2n_{k}}}\varphi \left( q_{n_{k}}\right) }.
\end{equation*}%
Hence we can write:

\begin{eqnarray}
&&\mu \left\{ x\in G_{m}:\frac{\left\vert \sigma
_{q_{_{_{n_{k}}}}}f_{_{_{_{n_{k}}}}}\left( x\right) \right\vert }{\varphi
\left( q_{_{_{n_{k}}}}\right) }\geq \frac{c}{M_{_{_{_{2n_{k}}}}}\varphi
\left( q_{_{_{n_{k}}}}\right) }\right\}  \label{16} \\
&\geq &\mu \left\{ x\in I_{_{2n_{k}}}^{2,4}:\frac{\left\vert \sigma
_{q_{_{_{n_{k}}}}}f_{_{_{_{n_{k}}}}}\left( x\right) \right\vert }{\varphi
\left( q_{_{_{n_{k}}}}\right) }\geq \frac{c}{M_{_{_{_{2n_{k}}}}}\varphi
\left( q_{_{_{n_{k}}}}\right) }\right\}  \notag \\
&\geq &\mu \left( I_{_{2n_{k}}}^{2,4}\right) >c>0.\qquad  \notag
\end{eqnarray}

From (\ref{15}) and (\ref{16}) we have

\begin{eqnarray*}
&&\frac{\frac{c}{M_{_{_{_{2n_{k}}}}}\varphi \left( q_{_{_{n_{k}}}}\right) }%
\left( \mu \left\{ x\in G_{m}:\frac{\left\vert \sigma
_{q_{_{_{n_{k}}}}}f_{_{_{_{n_{k}}}}}\left( x\right) \right\vert }{\varphi
\left( q_{_{_{n_{k}}}}\right) }\geq \frac{c}{M_{_{_{_{2n_{k}}}}}\varphi
\left( q_{_{_{n_{k}}}}\right) }\right\} \right) ^{1/p}}{\left\Vert
f_{n_{k}}\left( x\right) \right\Vert _{H_{p}}} \\
&\geq &\frac{c}{M_{_{_{_{2n_{k}}}}}\varphi \left( q_{_{_{n_{k}}}}\right)
M_{_{2n_{k}}}^{1-1/p}} \\
&=&c\frac{M_{_{2n_{k}}}^{1/p-2}}{\varphi \left( q_{_{_{n_{k}}}}\right) }%
\rightarrow \infty \text{,\qquad when }k\rightarrow \infty .
\end{eqnarray*}

Theorem 2 is proved.


\begin{thebibliography}{99}
\bibitem{AVD} G. N. AGAEV, N. Ya. VILENKIN, G. M. DZHAFARLY and A. I.
RUBINSHTEIN, Multiplicative systems of functions and harmonic analysis on
zero-dimensional groups, Baku, Ehim, 1981 (in Russian).

\bibitem{BGG} I. BLAHOTA, G. GÁT and U. GOGINAVA, Maximal operators of Fejer
means of double Vilenkin-Fourier series, Colloq. Math. 107 (2007), no. 2,
287--296.

\bibitem{BGG2} I. BLAHOTA, G. GÁT and U. GOGINAVA, Maximal operators of Fejé%
r means of Vilenkin-Fourier series. JIPAM. J. Inequal. Pure Appl. Math. 7
(2006), 1- 7 .

\bibitem{Fu} N. J. FUJII, A maximal inequality for $H^{1}$ functions on the
generalized Walsh-Paley group, Proc. Amer. Math. Soc. 77 (1979), lll-116.

\bibitem{gat} G. GÁT, Ces\`{a}ro means of integrable functions with respect
to unbounded Vilenkin systems. J. Approx. Theory 124 (2003), no. 1, 25--43

\bibitem{GoSzeged} U. GOGINAVA, Maximal operators of Fejér-Walsh means. Acta
Sci. Math. (Szeged) 74 (2008), no. 3-4, 615--624.

\bibitem{GoPubl} U. GOGINAVA, The maximal operator of the Fejér means of the
character system of the $p$-series field in the Kaczmarz rearrangement.
Publ. Math. Debrecen 71 (2007), no. 1-2, 43--55.

\bibitem{GoAMH} U. GOGINAVA, Maximal operators of Fejér means of double
Walsh-Fourier series. Acta Math. Hungar. 115 (2007), no. 4, 333--340

\bibitem{GNCz} U. GOGINAVA and K. NAGY , On the maximal operator of
Walsh-Kaczmarz-Fejer means, Czechoslovak Math. J. 62, 3 (2011), 673-686.

\bibitem{PS} J. PÁL and P. SIMON, On a generalization of the comncept of
derivate, Acta Math. Hung., 29 (1977), 155-164.

\bibitem{Sc} F. SCHIPP, Certain rearranngements of series in the Walsh
series, Mat. Zametki, 18 (1975), 193-201.

\bibitem{Si1} P. SIMON, Cesaro summability wish respect to two-parameter
Walsh sistems, Monatsh. Math .,131 (2000),321-334.

\bibitem{Si2} P. SIMON, Inverstigations wish respect to the Vilenkin sistem,
Annales Univ. Sci. Budapest Eotv., Sect. Math., 28 (1985) 87-101.

\bibitem{tep1} G. TEPHNADZE, Fejér means of Vilenkin-Fourier series, Studia
Scientiarum Mathematicarum Hungarica, (to appear).

\bibitem{tep2} G. TEPHNADZE, On the maximal operator of Vilenkin-Fejér
means, Turk. J. Math, (to appear).

\bibitem{Vi} N. Ya. VILENKIN, A class of complate ortonormal systems, Izv.
Akad. Nauk. U.S.S.R., Ser. Mat., 11 (1947), 363-400

\bibitem{We1} F. WEISZ, Martingale Hardy spaces and their applications in
Fourier Analysis, Springer, Berlin-Heideiberg-New York, 1994.

\bibitem{We2} F. WEISZ, Cesáro summability of one and two-dimensional
Fourier series, Anal. Math. 5 (1996), 353-367.

\bibitem{We3} F. WEISZ, Summability of multi-dimensional Fourier series and
Hardy space, Kluwer Academic, Dordrecht, 2002.

\bibitem{we4} F. WEISZ, Weak type inequalities for the Walsh and bounded
Ciesielski systems. Anal. Math. 30 (2004), no. 2, 147--160.

\bibitem{Zy} A. ZYGMUND, Trigonometric Series, Vol. 1, Cambridge Univ.
Press, 1959.
\end{thebibliography}
\end{document}